\documentclass{amsart}
\usepackage[english]{babel}
\usepackage{amssymb,textcomp,bbm,latexsym}
\usepackage{amsrefs}
\usepackage{amsthm}
\usepackage{color}
\newcommand{\assign}{:=}
\newcommand{\mathd}{\mathrm{d}}
\newcommand{\tmop}[1]{\ensuremath{\operatorname{#1}}}
\newenvironment{itemizedot}{\begin{itemize} }{\end{itemize}}

\theoremstyle{plain} 

\newtheorem{theorem}[equation]{Theorem}

\theoremstyle{definition}
\newtheorem{definition}[equation]{Definition} 

\theoremstyle{remark}

\numberwithin{equation}{section}

\begin{document}

\title{A matrix weighted bilinear Carleson Lemma and Maximal Function}

\author{Stefanie Petermichl}{\thanks{partially supported by ERC grant CHRiSHarMa no. DLV-682402}}

\author{Sandra Pott}

\author{Maria Carmen Reguera}

\

\begin{abstract}
  We prove a bilinear Carleson embedding theorem with matrix weight and
  scalar measure. In the scalar case, this becomes exactly the well known
  weighted bilinear Carleson embedding theorem. Although only allowing scalar Carleson
  measures, it is to date the only extension to the bilinear setting of the
  recent Carleson embedding theorem by Culiuc and Treil that features a matrix 
  Carleson measure and a matrix weight.
  
  It is well known that a Carleson embedding theorem implies a Doob's 
 maximal inequality and this holds true in the matrix weighted setting with 
 an appropriately defined maximal operator. It is also known that a dimensional 
 growth must occur in the Carleson embedding theorem with matrix Carleson measure, 
 even with trivial weight. We give a definition of a maximal type function whose 
 norm in the matrix weighted setting does not grow with dimension.
\end{abstract}

{\maketitle}

\

\section{Introduction}

The Carleson embedding theorem, first developed in the context of L. Carleson's celebrated proof of the Corona Theorem, is a classical theorem in harmonic analysis
with many applications to PDE. It states that a Carleson measure $\mu$ gives rise
to a bounded embedding $L^2(\mathbb{T}) \rightarrow L^2(\mathbb{D}, \mu)$. In other words, it states that $L^2$ boundedness holds
if and only if it holds when testing on certain elementary functions. 

In this note we are concerned with matrix analogues. An unweighted Carleson
lemma with matrix Carleson measure holds trivially, derived from the scalar
case. However, the dimensional growth was under investigation in the 90s and
it was observed that at least a logarithmic growth with dimension,  i.e. the size of the
arising matrices, must occur. The sharp dimensional estimate was found in
work of Nazarov, Pisier, Treil and Volberg \cite{MR1880830} by means of a clever trick with Bellman functions,
using elementary complex analysis (Zweikonstantensatz).

Recently Culiuc and Treil \cite{CT} added an observation to this idea to obtain
the matrix weighted version for an embedding theorem with matrix 
Carleson measures,  without any restrictions on the weight other than it being a matrix weight. Previous versions assumed the matrix $A_2$ condition on the weight, a strong
assumption, which did not appear in the scalar case.

The result of Culiuc and Treil implies a Doob's inequality in the matrix weighted setting. Since even the unweighted Carleson embedding theorem with matrix measure induces a dimensional growth with the size of the matrices, this growth is inherited by the thus obtained Doob's inequality for the maximal operator defined in Definition \ref{usualmax}. We are interested in a version of a maximal type operator that does not have growth with dimension, see Definition \ref{badmax} and  Theorem \ref{t.badmax}. This operator has certain nice properties, but it is not a classical maximal operator. The estimate we provide is via Bellman function.

A weighted bilinear version of the Carleson embedding lemma, taking two functions as 
arguments, was an important tool in the early days of sharp weighted theory or two--weighted theory of singular
operators in the scalar case. The bilinear Carleson embedding theorem (BET) states that a Carleson
measure gives rise to a bilinear estimate. \

Motivated by the so--called matrix $A_2$ conjecture, which asks for the exact
growth 
of the matrix weighted norm of the Hilbert transform acting on vector
functions in terms of the matrix $A_2$-constant , we aim for a bilinear
version of the matrix weighted Carleson lemma.  Boundedness of matrix weighted Hilbert transform was first studied by Treil and Volberg in \cite{MR1428818}.
We do not succeed entirely, in the sense
that we have to assume that our Carleson measure has scalar entries, see Theorem \ref{t.main}.
The restriction to scalar Carleson sequence  is serious, but it is not very surprising, given the array of difficulties
encountered in the task of finding various sharp weighted estimates in the
matrix case. As of today, the only optimal estimate for a singular operator is
for the square function with matrix weight by Hyt\"onen, Petermichl and Volberg \cite{HPV}
and Treil \cite{T}, and a matrix-weighted maximal function \cite{IKP}. The estimate for the Hilbert transform is still open, with the best
to date estimate  missing the
sharp conjecture by a half power of the characteristic, see Nazarov, Petermichl, Treil and Volberg \cite{MR3689742}. Despite the many
advances made in the scalar weighted theory, matrix weights do not seem to
respond well to most of the tools known to us today.

The first proof of a version of the scalar bilinear Carleson lemma is found
in \cite{MR1897034} and was rather complicated, using tools and ideas from
\cite{NTV}, featuring a Bellman function and three conditions on the
measure sequence rather than one. It was understood for some time by the
experts that two of the arising conditions were actually redundant. Two of the conditions stated in \cite{MR1897034} are implied by the third, as can be
seen by a simple Bellman function argument. Indeed, we prove that this implication
remains correct in the case of matrix weights and scalar Carleson measure, see Theorem \ref{t.redundant}.
Despite this, we believe that the original proof of the bilinear Carleson lemma
still presents obstacles in the matrix case, even with scalar coefficients,
such as considered here. The argument we apply in this note has the flavour of
level sets found in \cite{MR3406523}. 

We hope that the result on BET can be improved - we give some indications of
possible questions at the end of this note. Either improvement could be used
towards improving the bound on the Hilbert transform with matrix weight. Even
though Carleson type theorems with matrix weight and scalar measure do appear
to be useful in getting estimates for singular operators, for example in
\cite{MR1428818}, the constants increase by at least a logarithmic factor. In addition,
we do not know if our version of BET can be used directly for such estimates.
It seems to us at least some improvement would be required, see the end of
this note.

We also point out open questions arising on the dimensional dependence at the end of this note. 

\section{Notation and results}

Let us say $Q_0 = [0, 1]$ is endowed with a dyadic filtration generated by dyadic
intervals $\mathcal{D}=\mathcal{D} (Q_0)$. We index the time by non--negative integers $k\ge 0$ so
that $\mathcal{F}_0$ has the atom $\mathcal{D}_0=\{Q_0\}$ of size $| Q_0 | = 1^0$ and using atoms $\mathcal{D}_k =\{Q \in \mathcal{D} (Q_0) : | Q | = 2^{- k}\}$ we define the generated $\sigma$--algebra at time $k$: $\mathcal{F}_k =
\sigma (\mathcal{D}_k)$. 

We call a $d \times d$ matrix--valued function
$W$ a weight if $W$ is entry--wise locally integrable and if $W (x)$ is
positive semidefinite almost everywhere. One defines $L^2 (\mathbbm{C}^d;W)$
to be the set of vector functions with
$$ \| f \|^2_{L^2 (\mathbbm{C}^d;W)} = \int_{Q_0} \| W^{1 / 2} (x) f (x)
   \|_{\mathbbm{C}^d}^2 \mathd x = \int_{Q_0} \langle W (x) f (x), f (x)
   \rangle_{\mathbbm{C}^d} \mathd x < \infty .$$
Let us denote by $\langle \cdot \rangle_Q$ the average of a scalar, vector or
matrix function over $Q$ and let us denote by $\cdot(Q)$ the integral over $Q$. 

The precise statement of the bilinear matrix weighted Carleson lemma we prove is the following.

\begin{theorem}[BET] \label{t.main}
Let $(\alpha_Q)$ be a sequence of non-negative scalars. Then for vector functions  $f, g$
  supported in $Q_0$,
\begin{eqnarray*}
&& \frac{1}{| K |} \sum_{Q \in \mathcal{D} (K)} \alpha_Q \le A \quad \forall
     K \in \mathcal{D} (Q_0)\\
     & \Longrightarrow& \exists B=B(A): 
 \sum_{Q \in \mathcal{D} (Q_0)} \alpha_Q | \langle \langle W
     \rangle^{- 1}_Q \langle W^{1 / 2} f \rangle_Q, \langle W^{- 1} \rangle^{-
     1}_Q \langle W^{- 1 / 2} g \rangle_Q \rangle_{\mathbbm{C}^d} | \\
   && \quad \le B    \| f \|_{L^2 (\mathbbm{C}^d)}  \| g \|_{L^2 (\mathbbm{C}^d)} \quad \forall f,g. 
\end{eqnarray*}     
\end{theorem}

We recall the recent result \cite{CT} where the Carleson sequence is allowed to be matrix valued, but only one function $f$ is being used:
\begin{theorem}[CET]\label{theoremCT}
    (Culiuc--Treil) Let $W$ be an invertible selfadjoint matrix weight of size
    $d$ and $(A_Q)$ a sequence of positive semidefinite matrices of size $d$.
    Then for a vector function $f$ supported in $Q_0$,
 \begin{eqnarray*}
 && \frac{1}{| K |} \sum_{Q \in \mathcal{D} (K)} \langle W \rangle_Q A_Q
       \langle W \rangle_Q \le A \langle W \rangle_K \quad  \forall K \in
       \mathcal{D} (Q_0)\\
       &\Longrightarrow& \exists B=B(A) : \sum_{Q \in \mathcal{D} (Q_0)} \| A_Q^{1
       / 2} \langle W^{1 / 2} f \rangle_Q \|_{\mathbbm{C}^d}^2\\
       &&\quad \le  B  \| f \|^2_{L^2 (\mathbbm{C}^d)}  \quad \forall f  .
 \end{eqnarray*}      
  \end{theorem}
  
Related to Theorem \ref{t.main} is the following implication, which we prove via Bellman functions.

\begin{theorem}\label{t.redundant}
  Let $(\alpha_Q)$ be a non--negative sequence and $W$ be a matrix weight such
  that $W$ and $W^{- 1}$ are summable on all dyadic intervals. Then
$$ \frac{1}{| K |} \sum_{Q \in \mathcal{D} (K)} \alpha_Q \leqslant 1 \quad \forall
     K \in \mathcal{D} (Q_0)$$
		$$\Longrightarrow \frac{1}{| K |} \sum_{Q \in
     \mathcal{D} (K)} \alpha_Q \langle W^{- 1} \rangle_Q^{- 1} \leqslant 4
     \langle W \rangle_K \quad \forall K \in \mathcal{D} (Q_0), $$
  where the second inequality is in the sense of operators.
\end{theorem}

Its motivation and historic meaning are discussed in Section \ref{s.BET}.

  Here is the usual definition of the maximal function with matrix measure, with the
supremum over the dyadic intervals:

\begin{definition}\label{usualmax}
  Let
$$ M^W f (x) = \sup_{Q \in \mathcal{D} (Q_0), x \in Q} \| W^{1 / 2} (x)
     \langle W \rangle^{- 1}_Q \langle W^{1 / 2} f \rangle_Q
     \|_{\mathbbm{C}^d} . $$
\end{definition}

Definition \ref{usualmax} is one way of extending the martingale maximal function to the matrix weighted setting, see Section \ref{s.max} for the motivation. Indeed, using Theorem \ref{theoremCT} one can show

\begin{theorem}\label{t.usualmax}
  The operator $M^W : L^2 (\mathbbm{C}^d) \rightarrow L^2 (\mathbbm{C})$ is bounded. 
\end{theorem}

In turn, this estimate is used to obtain Theorem \ref{t.main}. The implied constants thus obtained for the estimate in Theorem \ref{t.usualmax} are no better than the constants obtained in Theorem \ref{theoremCT}. Consider the following alternative below. See Section \ref{s.max} for the motivation of this object.  

\begin{definition}\label{badmax}
For time $k\ge 0$ and $x \in Q \in \mathcal{D}_k$, let 
\[ M_k^W f (x) = \| \langle W \rangle^{1 / 2}_Q \langle W \rangle^{- 1}_{K_k
   (Q)} \langle W^{1/2} f \rangle_{K_k (Q)} \|_{\mathbbm{C}^d}  \]
with $K_k:\mathcal{D}_k \to \cup_{i=0}^k \mathcal{D}_i ;Q \mapsto K_k(Q)$ defined inductively: let   
$ K_0 (Q_0) = Q_0 $ and for $k>0$ let 
$$ K_k (Q) =  Q  \text{ if }  \| \langle W \rangle^{-1 / 2}_Q
     \langle W^{1/2} f \rangle_Q\|_{\mathbbm{C}^d} > \| \langle W \rangle^{1 / 2}_Q \langle W
     \rangle^{- 1}_{K_{k-1} (\hat{Q})} \langle W^{1/2} f \rangle_{K_{k-1} (\hat{Q})} \|_{\mathbbm{C}^d}$$
and $   K_k(Q)=K_{k-1} (\hat{Q})$ else. Here, $\hat{Q}$ denotes the parent of $Q$.
\end{definition}   
Here is our result, which carries the same norm estimate as seen in Doob's inequality and in particular does not depend on dimension.
\begin{theorem}\label{t.badmax}
  For all $k \geqslant 0$ there holds $\| M^W_k \|_{L^2 (\mathbbm{C}^d)
  \rightarrow L^2(\mathbb{C})} \leqslant 2.$
\end{theorem}

\section{Maximal function}\label{s.max}

In this section we prove that the matrix weighted maximal function is bounded without
additional assumptions on the matrix weight. This line of argument is well
known to the experts in the area and rests on the strength of the matrix
weighted Carleson lemma by Culiuc and Treil \cite{CT}. This estimate inherits the dimensional 
growth that occurs in the embedding theorem. We then proceed with the definition of a maximal function via an adapted sequence. We show that for this definition, the dimensional growth does not occur. The arguments rests on the Bellman function argument in \cite{NT}.

\subsection{Usual Maximal function}

Definition \ref{usualmax} is motivated as follows. Since maximal norm quantities have scalar values, the matrix weight is included in the definition. Let the expression
\[ (M_w f) (x) = \sup_{Q : x \in Q} | \langle f \rangle_{Q, w} | = \sup_{Q : x
   \in Q} | w (Q)^{-1} \int_Q f w | = \sup_{Q : x
   \in Q} |  \langle w \rangle_Q^{-1} \langle fw \rangle_Q  | \]
denote the classical dyadic maximal function with measure $w$ for non--negative
weight $w$ and scalar valued $f$. Observe that since $w$ is non--negative, the suprema
 in $M_w f (x)$ and $M^w f (x)$ are attained at the same interval $I \ni x$ and 
\[ \| M^w f (x) \|_{L^2} = \| M_w f (x) \|_{L^2 (w)} . \]
As a consequence of Doob's theorem, we have the estimate $$\| M^w f (x) \|_{L^2}=\| M_w f \|_{L^2
(w)} \leqslant 2 \| f \|_{L^2 (w)}=2 \| fw^{1/2} \|_{L^2}.$$ We turn to the proof of Theorem \ref{t.usualmax} 
via linearization, a well known argument. 

\begin{proof}
  For any $Q \in \mathcal{D} (Q_0)$ define $E_Q$ to be the subset of $x \in Q$
  so that $Q$ is the maximal cube with respect to the collection that attains
  the maximum in the definition of $M^W$. So
$$M^W f (x) = \sum_{Q \in \mathcal{D} (Q_0)} \| W^{1 / 2} (x) \langle W
     \rangle^{- 1}_Q \langle W^{1 / 2} f \rangle_Q \|_{\mathbbm{C}^d}
     \chi_{E_Q} (x) . $$
  Now by the disjointness of the $E_Q$
  \begin{eqnarray*}
    \| M^W f \|^2_{L^2 (\mathbbm{C})} & = & \int_{Q_0} \sum_{Q \in \mathcal{D}
    (Q_0)} \| W^{1 / 2} (x) \langle W \rangle^{- 1}_Q \langle W^{1 / 2} f
    \rangle_Q \|_{\mathbbm{C}^d}^2 \chi_{E_Q} (x) \mathd x\\
    & = & \sum_{Q \in \mathcal{D} (Q_0)} \int_{E_Q} \langle W (x) \langle W
    \rangle^{- 1}_Q \langle W^{1 / 2} f \rangle_Q, \langle W \rangle^{- 1}_Q
    \langle W^{1 / 2} f \rangle_Q \rangle_{\mathbbm{C}^d} \mathd x\\
    & = & \sum_{Q \in \mathcal{D} (Q_0)} \langle \langle W \rangle_{E_Q} |
    E_Q | \langle W \rangle^{- 1}_Q \langle W^{1 / 2} f \rangle_Q, \langle W
    \rangle^{- 1}_Q \langle W^{1 / 2} f \rangle_Q \rangle_{\mathbbm{C}^d} .
  \end{eqnarray*}
  To proceed with the estimate, we first observe that
$$ \frac{1}{| K |} \sum_{Q \in \mathcal{D} (K)} \langle W \rangle_Q \langle
     W \rangle^{- 1}_Q \langle W \rangle_{E_Q} | E_Q | \langle W \rangle^{-
     1}_Q \langle W \rangle_Q \leqslant \langle W \rangle_K $$
  in the sense of operators, which is satisfied thanks to the disjointness of
  the $E_Q$. Now we use Theorem \ref{theoremCT} with
$$ A_Q = \langle W \rangle^{- 1}_Q \langle W \rangle_{E_Q} | E_Q | \langle W
     \rangle^{- 1}_Q .$$
  The conclusion of this theorem gives us the desired estimate:
  \begin{eqnarray*}
    \| M^W f \|^2_{L^2 (\mathbbm{C})} & = & \sum_{Q \in \mathcal{D} (Q_0)}
    \langle \langle W \rangle^{- 1}_Q \langle W \rangle_{E_Q} | E_Q | \langle
    W \rangle^{- 1}_Q \langle W^{1 / 2} f \rangle_Q, \langle W^{1 / 2} f
    \rangle_Q \rangle_{\mathbbm{C}^d}\\
    & = & \sum_{Q \in \mathcal{D} (Q_0)} \langle A_Q \langle W^{1 / 2} f
    \rangle_Q, \langle W^{1 / 2} f \rangle_Q \rangle_{\mathbbm{C}^d}\\
    & \le  & B(1) \| f \|^2_{L^2 (\mathbbm{C}^d)} .
  \end{eqnarray*}
\end{proof}
Observe that the implied constant is the same as in the conclusion of Theorem \ref{theoremCT}.

\subsection{Maximal function with poor memory}

Let us now motivate and explain the different maximal function from Definition \ref{badmax}. 

Notice that the classical scalar maximal function is an adapted process if the
supremum is restricted to early times: adapted means that the expression
\[ (M_w f)_k (x) = \sup_{Q : x \in Q, | Q | \geqslant 2^{- k}} | \langle f
   \rangle_{Q, w} | \]
has the property $(M_w f)_k \in \mathcal{F}_k$. This notation means the function $(M_w f)_k$ is 
measurable in the filtration $\mathcal{F}_k$ for all $k$, which is so since it is piecewise constant on atoms $\mathcal{D}_k$. If we denote by $K_k (x)$ the
interval containing $x$ so that
\[ (M_w f)_k (x) = | \langle f \rangle_{K_k (x), w} |, \]
then \ $(M_w f)_{k + 1} (x) = \max \{ | \langle f \rangle_{K_k (x), w} |, |
\langle f \rangle_{Q, w} | \}$ where $| Q | = 2^{- k - 1}$ and $x \in Q$. It
therefore suffices to `memorize' just one past interval, so `poor memory' is
enough in the scalar case. These properties motivated our definition for the
matrix weighted case.

The definition above for the matrix weighted maximal operator $M^W$ is not
adapted when truncating because of the pointwise multiplication by $W (x)^{1 /
2}$. We now review the inductive definition of the adapted $M_k^W$ for
$k \geqslant 0$ with `poor memory' in this setting. In this case `poor memory'
may be `bad memory'.

For time $k\ge 0$ and $x \in Q \in \mathcal{D}_k$, the process we defined is  
\[ M_k^W f (x) = \| \langle W \rangle^{1 / 2}_Q \langle W \rangle^{- 1}_{K_k
   (Q)} \langle W^{1/2} f \rangle_{K_k (Q)} \|_{\mathbbm{C}^d}  \]
for dyadic intervals $K_k(Q)$ defined inductively with $ K_0 (Q_0) = Q_0 $
so that
\[ M_0^W f (x) = \| \langle W \rangle^{- 1 / 2}_{Q_0} \langle W^{1/2} f
   \rangle_{Q_0} \|_{\mathbbm{C}^d} = \| \langle W \rangle^{1 / 2}_{Q_0} \langle W \rangle^{-
   1}_{Q_0} \langle W^{1/2} f \rangle_{Q_0} \|_{\mathbbm{C}^d} . \]
Assuming $k \geqslant 1$ and $K_{k-1} (Q)$ chosen for all dyadic intervals $Q \in  \mathcal{D}_{k-1}$,
choose $ K_k (Q) \in \{ K_{k-1} (\hat{Q}), Q \}  $ for $Q \in \mathcal{D}_k$
by
$$ K_k (Q) =  Q  \text{ if }  \| \langle W \rangle^{-1 / 2}_Q 
     \langle W^{1/2} f \rangle_Q \|_{\mathbbm{C}^d} > \| \langle W \rangle^{1 / 2}_Q \langle W
     \rangle^{- 1}_{K_{k-1} (\hat{Q})} \langle W^{1/2} f \rangle_{K_{k-1} (\hat{Q})} \|_{\mathbbm{C}^d}$$
and $    K_k(Q)=K_{k-1} (\hat{Q}) $ otherwise.  Hence
\begin{eqnarray*}
&&M^W_k f(x) = \| \langle W \rangle^{1 / 2}_Q \langle W \rangle^{- 1}_{K_k (Q)} \langle W^{1/2} f
   \rangle_{K_k (Q)} \|_{\mathbbm{C}^d} \\
   &=& \max \left\{ \| \langle W \rangle^{1 / 2}_Q \langle W
   \rangle^{- 1}_{K_{k-1} (\hat{Q})} \langle W^{1/2} f \rangle_{K_{k-1} (\hat{Q})} \|_{\mathbbm{C}^d}, \|
   \langle W \rangle^{-1 / 2}_Q \langle W^{1/2} f \rangle_Q
   \|_{\mathbbm{C}^d} \right\} . 
\end{eqnarray*}

Observe that the time index $k$ on the function $K$ can be omitted without confusion. 
Notice that the sequence $(M_k^W f)_{k\ge 0}$ is adapted. As motivated by the scalar case, we
compete the value using the interval retained from the previous step against
the one from the finest filtration corresponding to $k$.

Let us discuss differences between $M_{\infty}^W f$ and $M^W$, say in the
case when $W$ and $f$ are dyadic step functions and measurable in some
$\mathcal{F}_k$. We observe that generally $M_{\infty}^W f$ is larger. If $Q$
is an atom in $\mathcal{F}_k$ and $x \in Q$ then $W (x) = \langle W \rangle_Q$
and for $x \in Q$ we have
\begin{eqnarray*}
 M_{\infty}^W f (x)  &=& \sup_{J : x \in J} \| \langle W \rangle_Q^{1 / 2}
   \langle W \rangle^{- 1}_J \langle W^{1/2} f \rangle_J \|_{\mathbbm{C}^d}  \\
  & \geqslant & \sup_{J : Q
   \subseteq J} \| \langle W \rangle_Q^{1 / 2} \langle W \rangle^{- 1}_J
   \langle W^{1/2} f \rangle_J \|_{\mathbbm{C}^d}  \\
   &\geqslant & M_k^W f (x) . 
 \end{eqnarray*}  
Indeed, the supremum taken in $M_k^W f (x)$ is a maximum over only two
competitors. To illustrate, let $x \in Q \subset \hat{Q}$ with $Q\in \mathcal{D}_k$. In the passage from $k-1$ to $k$, the first factor changes from $\langle W \rangle^{1 / 2}_{\hat{Q}}$ to $\langle W \rangle^{1 / 2}_Q$
which may mean that $K_{k-1} (\hat{Q})$ is not a sensible choice for a competitor at
all if one tried to estimate a `true' maximal function at time $k$. However in the scalar case, one obtains an approximating sequence for $M^wf$.

We proceed with the proof of Theorem \ref{t.badmax}.

\begin{proof}
  We modify the Bellman function argument from \cite{NT}. It shows a dimensionless bound
  for $M_k^W$ uniformly in $k$. Let
  \[ B (F, f, L, W) = 4 F - 2 \langle W^{- 1 / 2} f, L \rangle
     - 2 \langle L, W^{- 1 / 2} f \rangle + 2 \langle L, L \rangle .\] 
     We think of the variables having the following roles:
  \[ W_Q = \langle W \rangle_Q , f_Q = \langle W^{1/2} f \rangle_Q, F_Q =
     \langle \| f \|_{\mathbbm{C}^d} ^2 \rangle_Q \tmop{and} \]
  \[ L_Q = \langle W \rangle_Q^{1 / 2} \langle W \rangle^{- 1}_{K (Q)} \langle
     W^{1/2} f \rangle_{K (Q)}  \]
  with interval $K (Q) \in \{ K (\hat{Q}), Q \}$ inductively chosen for maximizing the norm 
  $\| \langle W \rangle_Q^{1 / 2} \langle W \rangle^{- 1}_{K (Q)} \langle W^{1/2} f
  \rangle_{K (Q)} \|_{\mathbbm{C}^d} $ as described above. 
  
  \paragraph{Domain of $B$.}The domain is given by $W_Q \geqslant 0$ and 
  \begin{equation}\label{e.upper.max}
  F_Q \geqslant \| W^{- 1 / 2}_Q f_Q \|_{\mathbbm{C}^d} ^2.
  \end{equation}
  To see this, recall the fact $$\langle W_Q^{-1} f_Q,f_Q\rangle_{\mathbb{C}^d}=\sup_{e \neq 0}\frac{|\langle f_Q,e\rangle_{\mathbb{C}^d} |^2}{\langle   W_Q e,e \rangle_{\mathbb{C}^d}},$$ for example from \cite{CT}. Then observe that since $W$ is selfadjoint,
  $$|\langle\langle W^{1/2}f,e\rangle_{\mathbb{C}^d} \rangle_Q|^2\le \langle \|f\|_{\mathbb{C}^d}^2 \rangle_Q \langle \langle We,e\rangle_{\mathbb{C}^d}\rangle_Q = F_Q \langle W_Q e,e \rangle_{\mathbb{C}^d}.$$
  Since $Q$ itself is an admissible choice for $K (Q)\in \{K(\hat{Q});Q  \}$ we require furthermore that
  \begin{equation}\label{e.intervalchoice}
  \| L_Q \|_{\mathbbm{C}^d}  \geqslant \| W^{- 1 / 2}_Q f_Q \|_{\mathbbm{C}^d} .
  \end{equation}
  
  \paragraph{Range of $B$.} Notice that $L_{Q_0} = W^{- 1 / 2}_{Q_0} f_{Q_0}$, since
  $Q_0$ is the only competitor. Therefore we have
  \[ B (F_{Q_0}, f_{Q_0}, L_{Q_0}, W_{Q_0}) = 4 F_{Q_0} - 2 \| W^{- 1 /
     2}_{Q_0} f_{Q_0} \|_{\mathbbm{C}^d} ^2 \leqslant 4 F_{Q_0} . \]
  We have also
  \[ B (F_Q, f_Q, L_Q, W_Q) \geqslant \| L_Q \|_{\mathbbm{C}^d} ^2 \]
  since
  \begin{eqnarray*}
    &&B (F_Q, f_Q, L_Q, W_Q)  \\
    & \geqslant & 4 
    \langle W^{- 1 / 2}_Q f_Q, W^{- 1 / 2}_Q f_Q \rangle_{\mathbb{C}^d} 
    - 2 \langle W^{- 1 / 2}_Q f_Q, L_Q \rangle_{\mathbb{C}^d}\\
    &  & 
    \quad - 2 \langle L_Q, W^{- 1 / 2}_Q f_Q \rangle_{\mathbb{C}^d} 
    + 2\langle L_Q, L_Q \rangle_{\mathbb{C}^d}\\
    & = & 
    \langle 2 W^{- 1 / 2}_Q f_Q - L_Q, 2 W^{- 1 / 2}_Q f_Q - L_Q\rangle_{\mathbb{C}^d}+ \langle L_Q, L_Q \rangle_{\mathbb{C}^d}\\
    & \geqslant &  \langle L_Q, L_Q \rangle_{\mathbb{C}^d}.
  \end{eqnarray*}
  
  \paragraph{Dynamics of $B$.}
  If we denote by $Q_{\pm}$ the left and right halves of $Q$, then $W_{Q_{\pm}}$ and $F_{Q_{\pm}}$ and $f_{Q_{\pm}}$ are
  the corresponding averages so that these variables are of martingale type.
  But
  \[ L_{Q_{\pm}} = \left\{ \begin{array}{ccc}
       W^{1 / 2}_{Q_{\pm}} W^{- 1 / 2}_Q L_Q  
       & \tmop{if} & 
       \| W^{- 1 /2}_{Q_{\pm}} f_{Q_{\pm}} \|_{\mathbbm{C}^d}  \leqslant \| W^{1 / 2}_{Q_{\pm}} W^{- 1 /
       2}_Q L_Q \|_{\mathbbm{C}^d}    \\
       W^{- 1 / 2}_{Q_{\pm}} f_{Q_{\pm}} & \tmop{if} & \| W^{- 1 /
       2}_{Q_{\pm}} f_{Q_{\pm}} \|_{\mathbbm{C}^d}  > \| W^{1 / 2}_{Q_{\pm}} W^{- 1 / 2}_Q L_Q
       \|_{\mathbbm{C}^d} 
     \end{array}\right.  . \]
  We claim that
  \[ B (F_Q, f_Q, L_Q, W_Q) \geqslant \frac{1}{2} B (F_{Q_+}, f_{Q_+},
     L_{Q_+}, W_{Q_+}) + \frac{1}{2} B (F_{Q_-}, f_{Q_-}, L_{Q_-}, W_{Q_-}) .
  \]
  
  \subparagraph{Case 1.}
  No changes in $K (Q) \rightarrow K (Q_{\pm})$, thus $K(Q_{\pm}) = K (Q)$. This implies 
  \begin{equation}\label{e.nochange}
    \| W^{- 1 /2}_{Q_{\pm}} f_{Q_{\pm}} \|_{\mathbbm{C}^d}  \leqslant \| W^{1 / 2}_{Q_{\pm}} W^{- 1 /
       2}_Q L_Q \|_{\mathbbm{C}^d} 
  \end{equation}
 and
  \[ L_{Q_{\pm}} = W^{1 / 2}_{Q_{\pm}} W^{- 1 / 2}_Q L_Q . \]
  Observe
  \begin{eqnarray*} 
  && \langle W_{Q_+} W^{- 1 / 2}_Q L_Q, W^{- 1 / 2}_Q L_Q \rangle_{\mathbbm{C}^d} 
  + \langle W_{Q_-} W^{- 1 / 2}_Q L_Q, W^{- 1 / 2}_Q L_Q \rangle_{\mathbbm{C}^d}  \\
     &=& 2 
 \langle W_Q W^{-1 / 2}_Q L_Q, W^{- 1 / 2}_Q L_Q \rangle_{\mathbbm{C}^d} 
 = 2 \langle L_Q, L_Q \rangle_{\mathbbm{C}^d} 
  \end{eqnarray*}   
  so
  \begin{eqnarray*}
    &  & 
    B (F_Q, f_Q, L_Q, W_Q)\\
    &  =& 
    4 F_Q - 2 \langle W^{- 1 / 2}_Q f_Q, L_Q \rangle_{\mathbbm{C}^d} 
    - 2 \langle L_Q, W^{- 1 / 2}_Q f_Q \rangle_{\mathbbm{C}^d} 
    + 2 \langle L_Q, L_Q \rangle_{\mathbbm{C}^d}\\
    & = & 
    4 F_Q - \langle f_{Q_+}, W^{- 1 / 2}_Q L_Q \rangle_{\mathbbm{C}^d} 
    - \langle f_{Q_-}, W^{- 1 / 2}_Q L_Q \rangle_{\mathbbm{C}^d}\\
    &  & 
    \quad - \langle W^{- 1 / 2}_Q L_Q, f_{Q_+}\rangle_{\mathbbm{C}^d}
    - \langle W^{- 1 / 2}_Q L_Q, f_{Q_-} \rangle_{\mathbbm{C}^d}\\
    &  & 
    \quad + \langle W_{Q_+} W^{- 1 / 2}_Q L_Q, W^{- 1 / 2}_Q L_Q \rangle_{\mathbbm{C}^d} 
    + \langle W_{Q_-} W^{- 1 / 2}_Q L_Q, W^{- 1 / 2}_Q L_Q \rangle_{\mathbbm{C}^d}\\
    & = & 
    2 F_{Q_+} - \langle W^{- 1 / 2}_{Q_+} f_{Q_+}, L_{Q_+} \rangle_{\mathbbm{C}^d} 
    - \langle L_{Q_+}, W^{- 1 / 2}_{Q_+} f_{Q_+} \rangle_{\mathbbm{C}^d} 
    + \langle L_{Q_+}, L_{Q_+} \rangle_{\mathbbm{C}^d}\\
    &  & 
    \quad +2 F_{Q_-} - \langle W^{- 1 / 2}_{Q_-} f_{Q_-}, L_{Q_-} \rangle_{\mathbbm{C}^d} 
    - \langle L_{Q_-}, W^{- 1 / 2}_{Q_-} f_{Q_-} \rangle_{\mathbbm{C}^d} 
    + \langle L_{Q_-}, L_{Q_-} \rangle_{\mathbbm{C}^d}\\
    & = & 
    \frac{1}{2} B (F_{Q_+}, f_{Q_+}, L_{Q_+}, W_{Q_+}) 
    + \frac{1}{2} B(F_{Q_-}, f_{Q_-}, L_{Q_-}, W_{Q_-}).
  \end{eqnarray*}
  We notice that in this case property (\ref{e.nochange}) was not used. 
  
  \subparagraph{Case 2.}
  Two changes in $K(Q) \rightarrow K(Q_{\pm})$ so $K
  (Q_{\pm}) = Q_{\pm}$. Let us thus assume
  \[ L_{Q_{\pm}} = W^{- 1 / 2}_{Q_{\pm}} f_{Q_{\pm}} . \]
  Due to the definition of $L_{Q_{\pm}}$ as a maximum of norms, we know that
  \begin{equation}\label{e.twochange}
   \| W^{- 1 / 2}_{Q_{\pm}} f_{Q_{\pm}} \|_{\mathbbm{C}^d} ^2 
   \geqslant \| W^{1 /2}_{Q_{\pm}} W^{- 1 / 2}_Q L_Q \|_{\mathbbm{C}^d} ^2 
 \end{equation}  
  and therefore
  \begin{eqnarray}\label{e.1}
   \| L_{Q_+} \|_{\mathbbm{C}^d}^2 + \| L_{Q_-} \|_{\mathbbm{C}^d}^2 
   &\geqslant &
    \| W^{1 / 2}_{Q_+} W^{- 1 /2}_Q L_Q \|_{\mathbbm{C}^d}^2 
    + \| W^{1 / 2}_{Q_-} W^{- 1 / 2}_Q L_Q \|_{\mathbbm{C}^d}^2\\
     \nonumber&=& 2 \| L_Q \|_{\mathbbm{C}^d}^2
     \geqslant 2 \| W^{- 1 / 2}_Q f_Q \|_{\mathbbm{C}^d}^2 ,
 \end{eqnarray}
   where we used property (\ref{e.twochange})  and the domain condition (\ref{e.intervalchoice}).  
  Observe that
  \begin{eqnarray*}
    &  & 
    2 \langle W^{- 1 / 2}_Q f_Q, W^{- 1 / 2}_Q f_Q \rangle_{\mathbbm{C}^d} \\
    & &
    \quad - 2 \langle W^{- 1 / 2}_Q f_Q, L_Q \rangle_{\mathbbm{C}^d} 
    - 2 \langle L_Q, W^{- 1 / 2}_Q f_Q \rangle_{\mathbbm{C}^d}
    + 2 \langle L_Q, L_Q \rangle_{\mathbbm{C}^d}\\
    & = & 
    2 \langle  W^{- 1 / 2}_Q f_Q - L_Q, 2 W^{- 1 / 2}_Q f_Q - L_Q\rangle_{\mathbbm{C}^d} 
    \geqslant 0,
  \end{eqnarray*}
  so using property (\ref{e.1})
  \begin{eqnarray*}
   && 
   - 2 \langle W^{- 1 / 2}_Q f_Q, L_Q \rangle_{\mathbbm{C}^d} 
   - 2 \langle L_Q, W^{- 1 / 2}_Q f_Q \rangle_{\mathbbm{C}^d} 
   + 2 \langle L_Q, L_Q \rangle_{\mathbbm{C}^d} \\
    & \geqslant & 
    - 2 \| W^{- 1 /2}_Q f_Q \|_{\mathbbm{C}^d}^2
     \geqslant 
    - \langle L_{Q_+}, L_{Q_+} \rangle_{\mathbbm{C}^d} 
    - \langle L_{Q_-}, L_{Q_-} \rangle_{\mathbbm{C}^d} .
  \end{eqnarray*}
  We will now show the dynamics inequality for this case:
  \begin{eqnarray*}
    &  & B (F_Q, f_Q, L_Q, W_Q)\\
    & = & 4 F_Q - 2 \langle W^{- 1 / 2}_Q f_Q, L_Q \rangle_{\mathbbm{C}^d} 
    - 2 \langle L_Q, W^{- 1 / 2}_Q f_Q \rangle_{\mathbbm{C}^d} 
    + 2 \langle L_Q, L_Q \rangle_{\mathbbm{C}^d}\\
    & \geqslant & 
    2 F_{Q_+} + 2 F_{Q_-} 
    - \langle L_{Q_+}, L_{Q_+} \rangle_{\mathbbm{C}^d} 
    - \langle L_{Q_-}, L_{Q_-} \rangle_{\mathbbm{C}^d}\\
    & = & 
    2 F_{Q_+} - \langle W^{- 1 / 2}_{Q_+} f_{Q_+}, L_{Q_+} \rangle_{\mathbbm{C}^d} 
    - \langle L_{Q_+}, W^{- 1 / 2}_{Q_+} f_{Q_+} \rangle_{\mathbbm{C}^d} 
    + \langle L_{Q_+}, L_{Q_+} \rangle_{\mathbbm{C}^d}\\
    &  & 
    + 2 F_{Q_-} - \langle W^{- 1 / 2}_{Q_-} f_{Q_-}, L_{Q_-} \rangle_{\mathbbm{C}^d} 
    - \langle L_{Q_-}, W^{- 1 / 2}_{Q_-} f_{Q_-} \rangle_{\mathbbm{C}^d} 
    + \langle L_{Q_-}, L_{Q_-} \rangle_{\mathbbm{C}^d}\\
    & = & 
    \frac{1}{2} B (F_{Q_+}, f_{Q_+}, L_{Q_+}, W_{Q_+}) 
    + \frac{1}{2} B (F_{Q_-}, f_{Q_-}, L_{Q_-}, W_{Q_-}) .
  \end{eqnarray*}
  
  \subparagraph{Case 3.}
  Mixed case, one change in $K (Q) \rightarrow K(Q_{\pm})$ with $K (Q_+) = K (Q)$ 
  and $K (Q_-) = K (Q)$. The other mixed case is symmetric. Let us thus assume
  \[ 
    L_{Q_-}  = W^{- 1 / 2}_{Q_-} f_{Q_-} 
   \; \tmop{and} \;
    L_{Q_+} = W^{1 / 2}_{Q_+} W^{- 1 / 2}_Q L_Q . 
  \]
  This gives 
  \[ 
     \| W^{- 1 / 2}_{Q_-} f_{Q_-} \|_{\mathbbm{C}^d} 
     > \| W^{1 / 2}_{Q_-} W^{- 1 / 2}_Q L_Q \|_{\mathbbm{C}^d}
    \; \tmop{and} \;
     \| W^{- 1 / 2}_{Q_+} f_{Q_+} \|_{\mathbbm{C}^d} 
     \leqslant 
     \| W^{1 / 2}_{Q_+}W^{- 1 / 2}_Q L_Q \|_{\mathbbm{C}^d} . 
 \]
  First notice
  \begin{eqnarray*}
    0 
    & \leqslant & 
    \langle W^{- 1 / 2}_{Q_-} f_{Q_-} - W^{1 / 2}_{Q_-} W^{- 1/ 2}_Q L_Q, W^{- 1 / 2}_{Q_-} f_{Q_-} 
    - W^{1 / 2}_{Q_-} W^{- 1 / 2}_Q L_Q\rangle_{\mathbbm{C}^d}\\
    & = & 
    \langle W^{- 1 / 2}_{Q_-} f_{Q_-}, W^{- 1 / 2}_{Q_-} f_{Q_-} \rangle_{\mathbbm{C}^d} 
    + \langle W^{1 / 2}_{Q_-} W^{- 1 / 2}_Q L_Q, W^{1 / 2}_{Q_-} W^{-1 / 2}_Q L_Q \rangle_{\mathbbm{C}^d}\\
    &  & 
    \quad - \langle W^{- 1 / 2}_{Q_-} f_{Q_-}, W^{1 / 2}_{Q_-} W^{- 1 / 2}_Q L_Q \rangle_{\mathbbm{C}^d} 
    - \langle W^{1 / 2}_{Q_-} W^{- 1 / 2}_Q L_Q, W^{- 1 / 2}_{Q_-}f_{Q_-} \rangle_{\mathbbm{C}^d}\\
    & = & 
    \langle W^{- 1 / 2}_{Q_-} f_{Q_-}, W^{- 1 / 2}_{Q_-} f_{Q_-}\rangle_{\mathbbm{C}^d} 
    + \langle W_{Q_-} W^{- 1 / 2}_Q L_Q, W^{- 1 / 2}_Q L_Q \rangle_{\mathbbm{C}^d}\\
    &  & 
    \quad - \langle f_{Q_-}, W^{- 1 / 2}_Q L_Q \rangle_{\mathbbm{C}^d} 
    - \langle W^{- 1 / 2}_Q L_Q, f_{Q_-} \rangle_{\mathbbm{C}^d} 
    + 2 \langle L_Q, L_Q \rangle_{\mathbbm{C}^d} \\
    &  &  
    \quad - \langle W_{Q_+} W^{- 1 / 2}_Q L_Q, W^{- 1 / 2}_Q L_Q \rangle_{\mathbbm{C}^d} 
    - \langle W_{Q_-} W^{- 1 / 2}_Q L_Q, W^{- 1/2}_Q L_Q \rangle_{\mathbbm{C}^d} .
  \end{eqnarray*}
  Thus
  \begin{eqnarray*}
    2 \langle L_Q, L_Q \rangle_{\mathbbm{C}^d} 
    & \geqslant & 
    \langle W_{Q_+} W^{- 1 / 2}_Q L_Q, W^{- 1 / 2}_Q L_Q \rangle_{\mathbbm{C}^d} 
    - \langle W^{- 1 / 2}_{Q_-} f_{Q_-}, W^{- 1/2}_{Q_-} f_{Q_-} \rangle_{\mathbbm{C}^d}\\
    &  & 
    \quad + \langle f_{Q_-}, W^{- 1 / 2}_Q L_Q \rangle_{\mathbbm{C}^d} 
    + \langle W^{- 1 / 2}_Q L_Q, f_{Q_-} \rangle_{\mathbbm{C}^d}
  \end{eqnarray*}
  and therefore
  \begin{eqnarray*}
    &  & 
    B (F_Q, f_Q, L_Q, W_Q)\\
    & = & 
    4 F_Q - 2 \langle W^{- 1 / 2}_Q f_Q, L_Q \rangle_{\mathbbm{C}^d} - 
    2 \langle L_Q, W^{- 1 / 2}_Q f_Q \rangle_{\mathbbm{C}^d} 
    + 2 \langle L_Q, L_Q \rangle_{\mathbbm{C}^d}\\
    & \geqslant & 
    4 F_Q  - \langle f_{Q_+}, W^{-1 / 2}_Q L_Q \rangle_{\mathbbm{C}^d} 
    - \langle f_{Q_-}, W^{- 1 / 2}_Q L_Q \rangle_{\mathbbm{C}^d}\\
    &  & 
    \quad - \langle W^{- 1 / 2}_Q L_Q, f_{Q_+}\rangle_{\mathbbm{C}^d} 
    - \langle W^{- 1 / 2}_Q L_Q, f_{Q_-} \rangle_{\mathbbm{C}^d}\\
    &  & 
    \quad + \langle f_{Q_-}, W^{- 1 / 2}_Q L_Q \rangle_{\mathbbm{C}^d} 
    + \langle W^{- 1 / 2}_Q L_Q, f_{Q_-} \rangle_{\mathbbm{C}^d}\\
    &  & 
    \quad +  \langle W_{Q_+} W^{- 1 / 2}_Q L_Q, W^{- 1 / 2}_Q L_Q \rangle_{\mathbbm{C}^d}  
    - \langle W^{- 1 / 2}_{Q_-} f_{Q_-}, W^{- 1/2}_{Q_-} f_{Q_-} \rangle_{\mathbbm{C}^d}\\
    & = &  
    2 F_{Q_+} + 2 F_{Q_-}\\
    &  &  
   \quad  - \langle W^{- 1 / 2}_{Q_+} f_{Q_+}, W^{1 / 2}_{Q_+} W^{- 1 / 2}_Q L_Q \rangle_{\mathbbm{C}^d} 
    - \langle W^{1 / 2}_{Q_+} W^{- 1 / 2}_Q L_Q, W^{- 1 / 2}_{Q_+} f_{Q_+} \rangle_{\mathbbm{C}^d}\\
    &  & 
    \quad + \langle W_{Q_+} W^{- 1 / 2}_Q L_Q, W^{- 1 / 2}_Q L_Q \rangle_{\mathbbm{C}^d} 
    - \langle W^{-1/2}_{Q_-} f_{Q_-}, W^{- 1 / 2}_{Q_-} f_{Q_-} \rangle_{\mathbbm{C}^d}\\
    & = & 
    2 F_{Q_+} - \langle W^{- 1/2}_{Q_+} f_{Q_+}, L_{Q_+} \rangle_{\mathbbm{C}^d} 
    - \langle L_{Q_+}, W^{- 1 / 2}_{Q_+}f_{Q_+} \rangle_{\mathbbm{C}^d} 
    + \langle L_{Q_+}, L_{Q_+} \rangle_{\mathbbm{C}^d}\\
    &  & 
    \quad + 2 F_{Q_-}  - \langle W^{- 1 / 2}_{Q_-} f_{Q_-}, L_{Q_-} \rangle_{\mathbbm{C}^d} 
    - \langle L_{Q_-}, W^{- 1 / 2}_{Q_-} f_{Q_-} \rangle_{\mathbbm{C}^d} 
    + \langle L_{Q_-}, L_{Q_-} \rangle_{\mathbbm{C}^d}\\
    & = & 
    \frac{1}{2} B (F_{Q_+}, f_{Q_+}, L_{Q_+}, W_{Q_+}) 
    + \frac{1}{2} B(F_{Q_-}, f_{Q_-}, L_{Q_-}, W_{Q_-}) .
  \end{eqnarray*}
Using the dynamics and size properties we obtain the inequality estimating the norm of $M_k^W$ uniformly in $k$:
  \begin{eqnarray*}
   \sum_{J \in \mathcal{D}_k( Q_0 )} \| L_J \|_{\mathbbm{C}^d} ^2 | J | 
   &\leqslant & 
   \sum_{J \in \mathcal{D}_k ( Q_0 )} B (F_J, f_J, L_J, W_J) | J | \\
   & \leqslant & 
   B(F_{Q_0}, f_{Q_0}, L_{Q_0}, W_{Q_0}) | Q_0 | \\
   & \leqslant & 
   4 F_{Q_0} | Q_0 |. 
  \end{eqnarray*}
  
\end{proof}

\section{Embeddings} \label{s.BET}

In this section we turn to the BET and the motivation and proof of the reduction estimate Theorem \ref{t.redundant}.

\subsection{Bilinear Embedding Theorem (BET)}\label{ss.BET}

We turn to the proof of Theorem \ref{t.main}.
\begin{proof}
 Let ${\mu} (\mathcal{K}) = \sum_{Q \in \mathcal{K}} \alpha_Q$ for any
  collection $\mathcal{K}$ of dyadic cubes. Let $F$ be any nonnegative function defined on
  the dyadic cubes. Then let  $\{ F (Q) > \lambda \}$ denote the collection of
  cubes $Q$ such that $F (Q) > \lambda$. It follows that
 $$ \int^{\infty}_0 {\mu} (\{ F (Q) > \lambda \}) \mathd \lambda = \sum_{Q \in
     \mathcal{D} (Q_0)} F (Q) \alpha_Q ,$$
  which is the classical fact on Choquet integrals.
  
  Let us define
$$ F (Q) = |\langle \langle W \rangle^{- 1}_Q \langle W^{1 / 2} f \rangle_Q,
     \langle W^{- 1} \rangle^{- 1}_Q \langle W^{- 1 / 2} g \rangle_Q
     \rangle_{\mathbbm{C}^d} |$$
  and for any $\lambda > 0$, let $\mathcal{J}_{\lambda}$ denote the
  collection of maximal dyadic intervals for which $F (Q) > \lambda$. Hence
$$ \sum_{Q \in \mathcal{J}_{\lambda}} \alpha_Q \leqslant \sum_{Q \in
     \mathcal{J}_{\lambda}} \sum_{Q' \in \mathcal{D} (Q)} \alpha_{Q'} \leqslant
     \sum_{Q \in \mathcal{J}_{\lambda}} | Q | . $$
  Now let
$$ \Phi (x) = M^W f (x) M^{W^{- 1}} g (x) .$$
  Observe that with
$$ A^W_Q f (x) \assign W^{1 / 2} (x) \langle W \rangle^{- 1}_Q \langle W^{1
     / 2} f \rangle_Q $$
  and
$$ A^{W^{- 1}}_Q g (x) = W^{- 1 / 2} (x) \langle W^{- 1} \rangle^{- 1}_Q
     \langle W^{- 1 / 2} g \rangle_Q, $$
   Cauchy Schwarz yields for all $x$ and all $Q$
 $$ F (Q) \leqslant \| A^W_Q f (x) \|_{\mathbbm{C}^d} \| A^{W^{- 1}}_Q g (x)
     \|_{\mathbbm{C}^d} . $$
  Now if $x \in Q$ with $F (Q) > \lambda$, then $\Phi (x) > \lambda$. Hence
 $$ \sum_{Q \in \mathcal{J}_{\lambda}} | Q | \leqslant | \{ x \in
     \mathbbm{R}: \Phi (x) > \lambda \} | .$$
  Integrating with respect to $d \lambda$ gives
 $$ \sum_{Q \in \mathcal{D} (Q_0)} \alpha_Q F (Q) \leqslant \int^{\infty}_0 |
     \{ x \in \mathbbm{R}: \Phi (x) > \lambda \} | \mathd \lambda =
     \int_{Q_0} M^W f (x) M^{W^{- 1}} g (x) \mathd x $$
  Using the above estimate of the maximal function in Theorem \ref{t.usualmax} and an application of
  Cauchy Schwarz finishes the estimate.
\end{proof}

\subsection{Altered Carleson condition}\label{ss.redundant}

As mentioned earlier, our BET reduces to the classical version of the scalar
case. The first proof of the scalar theorem was by the Bellman method, using a
rather cleverly built function, that in a way allows to deduce the bilinear
version BET from the linear one CET, see \cite{NTV}, \cite{MR1897034} and \cite{MR2354322}. This is exactly what we do here, too,
although with completely different methods, that are applicable to the
matrix weighted case. The original scalar version of BET, first formulated in \cite{NTV},
seemed to require three Carleson conditions, namely
$$
\frac{1}{| K |} \sum_{Q \in \mathcal{D} (K)} \frac{\alpha_Q}{\langle w\rangle_{Q}} \leqslant C\langle w^{-1}\rangle_{K}\quad  \forall
     K \in \mathcal{D} (Q_0),
$$

$$
\frac{1}{| K |} \sum_{Q \in \mathcal{D} (K)} \frac{\alpha_Q}{\langle w^{-1}\rangle_{Q}} \leqslant C\langle w\rangle_{K} \quad \forall
     K \in \mathcal{D} (Q_0),
$$

$$
\frac{1}{| K |} \sum_{Q \in \mathcal{D} (K)} \alpha_Q \leqslant C \quad \forall
     K \in \mathcal{D} (Q_0). 
$$

As mentioned in the introduction, it
turns out that the first two conditions can be removed in the scalar case
through the use of an embedding--like Bellman function, thus only retaining
one Carleson condition, the analog of what we used here. This is implicit in
 \cite{DP} and similar considerations also appeared in  \cite{MR2433959}. 
 This reduction also holds in the case of a scalar Carleson
measure and a matrix weight, as stated in Theorem \ref{t.redundant}.  We turn to its proof.

\begin{proof}
  Consider the matrix valued Bellman function of positive matrix variables $U, V$ and
  scalar variable $m$
 $$ B (U, V, m) = U - (m + 1)^{- 1} V^{- 1} . $$ The variables we have in mind are matrix valued 
 $U_K = \langle W \rangle_K$, $V_K = \langle W^{- 1}
  \rangle_K$ and scalar valued $m_K = | K |^{- 1} \sum_{Q \in \mathcal{D} (K)} \alpha_Q$.
  \paragraph{Domain of $B$.} This function has the domain conditions 
  $$\tmop{Id} \leqslant V^{1 / 2} U V^{1 / 2} \text{ and } 0 \leqslant m \leqslant 1.$$ 
  Indeed, the first condition is equivalent to     $U \geqslant V^{- 1}$
  and is implied by Lemma
  3.2 in \cite{MR1428818}: if $W$ is a matrix weight such that $W$ and $W^{-
  1}$ are summable on $Q$, then for all vectors $e$, $$\langle
  \langle W \rangle_Q e, e \rangle \geqslant \langle  \langle W^{- 1}
  \rangle^{- 1}_Q e, e \rangle .$$
  \paragraph{Range of $B$}
  We have the size estimate
\begin{equation}\label{e.size.redundant}
 0 \leqslant B (U, V, m) \leqslant U. 
 \end{equation}
Indeed, $0 \leqslant U - V^{- 1} \leqslant U - (m + 1)^{- 1} V^{- 1}
  \leqslant U.$
  
  \paragraph{Dynamics of $B$.} The function $B$ is concave: 
  Dropping the linear dependence on $U$, its Hessian acting on the matrix difference $\Delta V$ and scalar $\Delta m$
  is
 $$ - 2 V^{- 1} \Delta V V^{- 1} \Delta V V^{- 1} (m + 1)^{- 1} - 2 V^{- 1}
     \Delta V V^{- 1} (m + 1)^{- 2} \Delta m - 2 (m + 1)^{- 3} V^{- 1} (\Delta
     m)^2 . $$
  Observe that
$$ V^{- 1} \Delta V V^{- 1} \Delta V V^{- 1} (m + 1)^{- 1} \geqslant 0, $$
   and 
$$(m + 1)^{- 3} V^{- 1} (\Delta m)^2 \geqslant 0. $$
  Now add these two positive terms, then factor $- (m + 1)^{- 1}$ and reverse the sign:
  \begin{eqnarray*}
    &&  V^{- 1} \Delta V V^{- 1} \Delta V V^{- 1} + 2 V^{- 1} \Delta V V^{-
    1} (m + 1)^{- 1} \Delta m + (m + 1)^{- 2} V^{- 1} (\Delta m)^2\\
     &=& (\Delta
    m (m + 1)^{- 1}\tmop{Id}+V^{- 1} \Delta V )V^{- 1} (\Delta V V^{- 1} + \tmop{Id} (m + 1)^{- 1}
    \Delta m)    
   \geqslant  0.
  \end{eqnarray*}
  This inequality proves concavity of $B$, which in turn implies mid--point concavity 
  \begin{eqnarray}\label{e.concavity.redundant}
  && B(U,V,m)\\
  \nonumber &\ge&\frac12 B(U+\Delta U,V+ \Delta U,m+ \Delta m) + \frac12 B(U-  \Delta U,V- \Delta V,m- \Delta m).
 \end{eqnarray} 
  
We have also
 \begin{equation}\label{e.der} 
( \partial B / \partial m )(U,V,m)= (m + 1)^{- 2} V^{- 1} \geqslant \frac14 V^{- 1}. 
 \end{equation}
  
  With martingale matrix variables $U_K = \langle W \rangle_K$, $V_K = \langle W^{- 1}
  \rangle_K$ and scalar variable $m_K = | K |^{- 1} \sum_{Q \in \mathcal{D} (K)} \alpha_Q$ in the
  domain, we see that 
  \begin{equation}\label{e.smartingale}
  m_K - | K |^{- 1} \alpha_K = \frac12 (m_{K_-} +
  m_{K_+}).
  \end{equation}
 
   The usual argument gives the estimate of the operator sum: Fix $K
  \in \mathcal{D} (Q_0)$ and estimate
  \begin{eqnarray*}
    | K | \langle W \rangle_K & = & | K | U_K\\
    & \geqslant & | K | B (U_K, V_K, m_K)\\
    & = & | K | B (U_K, V_K, m_K) - | K | B (U_K, V_K, m_K - | K |^{- 1} \alpha_K)\\
		& & \quad
    + | K | B (U_K, V_K, m_K - | K |^{- 1} \alpha_K)\\
    & \geqslant & \frac14 V_K^{- 1} \alpha_K  + | K | B (U_K, V_K, m_K - | K
    |^{- 1} \alpha_K)\\
    & \geqslant & \frac14 V_K^{- 1} \alpha_K + | K_- | B (U_{K_-}, V_{K_-},
    m_{K_-}) + | K_+ | B (U_{K_+}, V_{K_+}, m_{K_+}).
  \end{eqnarray*}
  In these inequalities we used the size estimate (\ref{e.size.redundant}), the intermediate value theorem together with the derivative estimate (\ref{e.der}) and the dynamics of the martingale variables and (\ref{e.smartingale}), and the 
  concavity of $B$ in the form of mid-point concavity equation (\ref{e.concavity.redundant}). Iterating this
  argument gives the desired estimate.
\end{proof}

\section{Questions}

\

Our theorems leave room for improvement in several directions.

\subsection{$A_2$ conjecture}

Currently, the matrix $A_2$ conjecture for the Hilbert transform is under investigation. This question asks for the growth of the function $\Phi$ so that 
$$ \| H f \|_{L^2(\mathbb{C}^d;W)}\le C \Phi([W]_{A_2})\| f\|_{L^2(\mathbb{C}^d;W)},$$
where the Hilbert transform $H$ is applied to the vector function $f$ componentwise and the matrix $A_2$ characteristic of Treil and Volberg \cite{MR1428818} is defined as
$$[W]_{A_2}=\sup_{Q} \| \langle W \rangle_Q^{1/2}  \langle W^{-1} \rangle_Q^{1/2}  \|^2.$$
In the scalar case the best growth is linear, and this is optimal, \cite{MR2354322}. In the original text on the matrix case \cite{MR1428818} that laid the ground for these investigations, the dependence on $[W]_{A_2}$ was not tracked. The first quantified result with better estimates is \cite{BPW}, which obtains the power $3/2$ with an additional logarithmic term. In the text \cite{MR3689742} the logarithmic term was dropped and the best known estimate of power $3/2$  stands at current time. To improve this power further, one may aim at some improvements in the embedding theorems we demonstrated in this note:
\begin{itemizedot}
  \item Even within the framework of scalar $\alpha_Q$ in BET it would be useful to be able to
  get an estimate of the (larger) sum
 $$ \sum_{Q \in \mathcal{D} (Q_0)} \alpha_Q \| \langle W \rangle^{- 1}_Q
     \langle W^{1 / 2} f \rangle_Q \|_{\mathbb{C}^d} \| \langle W^{- 1} \rangle^{- 1}_Q
     \langle W^{- 1 / 2} g \rangle_Q \|_{\mathbb{C}^d} $$
  as this may improve the matrix weighted estimate for certain classical
  operators. Staying with scalar coefficients is unlikely to press the
  constants to the desired linear estimate  in the matrix $A_2$ constant, but one may hope for just an extra
  logarithmic term instead of the extra half power in current estimates.
  
  \item Does Theorem \ref{t.redundant} hold with a matrix Carleson measure? From the required
  size conditions a natural choice for a Bellman function is $$U - V^{- 1 / 2}
  (M + \tmop{Id})^{- 1} V^{- 1 / 2}.$$ It is not clear to us if, for example,
  the convexity in $V$ is as required. Indeed, to obtain derivatives of square
  roots, one may use the spectral resolution formula.
  
  \item Find a version of BET with matrix coefficients $\alpha_Q$. 
 \end{itemizedot} 
  
  \subsection{Dimensional growth}
  
  In the 90s the dimensional growth of matrix versions of typical scalar results was under investigation. There was a series of results on the Carleson embedding theorem with matrix measure and the papaproducts or Hankel operators with matrix symbol \cite{K}, \cite{P}, \cite{MR1880830}. Some of the recent results indicate that dimensional growth may be a finer indicator of non--commutativity than the growth with the $A_2$ constant.

\begin{itemizedot}

  \item Find the dimensional growth of the matrix weighted CET Theorem \ref{theoremCT}. When the weight is the identity, then the dimensional growth is $\log ^2(d)$ and this is sharp, see \cite{MR1880830}. It is not clear if the upper estimate is attainable in the presence of a non--trivial matrix weight. The best to date estimate is of order $d^2$. It appears that the matrix weight induces additional non--commutativity that causes the clever argument in \cite{MR1880830} to fail.
  
 \item What is the dimensional growth of the matrix maximal function of Definition \ref{usualmax}? 
  
\end{itemizedot}


\begin{bibsection}
\begin{biblist}

\bib{MR2433959}{article}{
   author={Beznosova, O.},
   title={Linear bound for the dyadic paraproduct on weighted Lebesgue space
   $L_2(w)$},
   journal={J. Funct. Anal.},
   volume={255},
   date={2008},
   number={4},
   pages={994--1007},
}


 \bib{BPW}{article}{
 author={Bickel, K.},
 author={Petermichl, S.},
 author={Wick, B.},
  title={Bounds for the Hilbert transform with matrix $A_2$ weights},
  pages={1719--1743},
  date={2016},
  number={5},
  volume={250},
}


\bib{CT}{article}{
   author={Culiuc, A.},
   author={Treil, S.},
   title={The Carleson embedding theorem with matrix weights},
   journal={To appear in IMRN},
}

\bib{DP}{article}{
   author={Domelevo, K.},
   author={Petermichl, S.},
   title={Differential subordination under change of law},
   journal={To appear in Annals of Probability},
   
}


\bib{HPV}{article}{
   author={Hyt\"onen, T.},
   author={Petermichl, S.},
   author={Volberg, A.},
   title={The sharp square function estimate with matrix weight},
   journal={Submitted},
}

\bib{IKP}{article}{
author={J. Isralowitz}, 
author={H.-K. Kwon}, 
author={S.Pott},
title={Matrix-weighted norm inequalities for commutators and paraproducts with matrix symbols},
year={2017},
volume={96},
number={1},
 journal = {J. London Math. Soc.},
 pages={243--270}
}


\bib{K}{article}{
   author={Katz,N.},
   title={Matrix valued paraproducts},
   journal={J. of Fourier Anal. Appl.},
   date={1997},
   pages={913--921},
   volume={300},
}   


\bib{MR3689742}{article}{
   author={Nazarov, F.},
   author={Petermichl, S.},
   author={Treil, Sergei},
   author={Volberg, A.},
   title={Convex body domination and weighted estimates with matrix weights},
   journal={Adv. Math.},
   volume={318},
   date={2017},
   pages={279--306},
}

\bib{MR1880830}{article}{
   author={Nazarov, F.},
   author={Pisier, G.},
   author={Treil, S.},
   author={Volberg, A.},
   title={Sharp estimates in vector Carleson imbedding theorem and for
   vector paraproducts},
   journal={J. Reine Angew. Math.},
   volume={542},
   date={2002},
   pages={147--171},
}

\bib{NT}{article}{
   author={Nazarov, F.},
   author={Treil, S.},
   title={The hunt for a Bellman functions: applications to estimates of singular integral operators and to other classical problems in harmonic analysis},
   journal={Algebra i Analysis},
   volume={8},
   date={1997},
   number={5},
   pages={32--162},
}

	

\bib{NTV}{article}{
   author={Nazarov, F.},
   author={Treil, S.},
   author={Volberg, A.},
   title={The Bellman functions and two-weight inequalities for Haar
   multipliers},
   journal={J. Amer. Math. Soc.},
   volume={12},
   date={1999},
   number={4},
   pages={909--928},
}



\bib{NTVbis}{article}{
   author={Nazarov, F.},
   author={Treil, S.},
   author={Volberg, A.},
   title={Two weight inequalities for individual Haar multipliers and other
   well localized operators},
   journal={Math. Res. Lett.},
   volume={15},
   date={2008},
   number={3},
   pages={583--597},
}


\bib{P}{article}{
   author={Petermichl, S.},
   title={Dyadic Shifts and a Logarithmic Estimate for Hankel Operators with Matrix Symbol},
   date={2000},
   pages={455--460},
   number={1},
   volume={330},
   journal={Comptes Rendus Acad. Sci. Paris},
}   


\bib{MR2354322}{article}{
   author={Petermichl, S.},
   title={The sharp bound for the Hilbert transform on weighted Lebesgue
   spaces in terms of the classical $A_p$ characteristic},
   journal={Amer. J. Math.},
   volume={129},
   date={2007},
   number={5},
   pages={1355--1375},
}


\bib{MR1897034}{article}{
   author={Petermichl, S.},
   author={Wittwer, J.},
   title={A sharp estimate for the weighted Hilbert transform via Bellman
   functions},
   journal={Michigan Math. J.},
   volume={50},
   date={2002},
   number={1},
   pages={71--87},
}

\bib{MR3406523}{article}{
   author={Thiele, C.},
   author={Treil, S.},
   author={Volberg, A.},
   title={Weighted martingale multipliers in the non-homogeneous setting and
   outer measure spaces},
   journal={Adv. Math.},
   volume={285},
   date={2015},
   pages={1155--1188},
}

\bib{T}{article}{
   author={Treil, S.},
   title={Mixed $A_{2}-A_{\infty}$ estimates of the non-homogeneous vector square function with matrix weights},
   journal={Submitted},
}

\bib{MR1428818}{article}{
   author={Treil, S.},
   author={Volberg, A.},
   title={Wavelets and the angle between past and future},
   journal={J. Funct. Anal.},
   volume={143},
   date={1997},
   number={2},
   pages={269--308},
}
		

\end{biblist}
\end{bibsection}
\end{document}